\documentclass[12pt]{article}
\usepackage{latexsym}
\usepackage{amsmath,amssymb}
\usepackage{bm}
\usepackage{graphicx,xcolor}
\usepackage{amssymb}
\usepackage{cite}
\usepackage{amsfonts, amscd, amsthm}
\usepackage[all]{xy}
\usepackage{here}
\usepackage{ulem}
\usepackage{mathabx}
\usepackage[titletoc,title]{appendix}
\usepackage{mathbbol}
\usepackage{url}

\renewcommand{\theequation}{\thesection.\arabic{equation}}

\newcommand{\lam}{\lambda}
\newcommand{\MGamma}{\mathit\Gamma}

\makeatletter
\@addtoreset{equation}{section}
\renewcommand{\theequation}{\thesection.\@arabic\c@equation}
\makeatother

\makeatletter
\renewcommand\appendix{\par
  \setcounter{section}{0}%
  \setcounter{subsection}{0}%
  \gdef\thesection{Appendix \@Alph\c@section }
  \renewcommand{\theequation}
  {\Alph{section}.\arabic{equation}}
}
\makeatother

\makeatletter
\newcounter{subeqncnt}
\def\thesubeqncnt{\alph{subeqncnt}}
\def\subequations{\begingroup%
\stepcounter{equation}\edef\@tempa{\theequation}%
\let\c@equation\c@subeqncnt\c@subeqncnt\z@
\edef\theequation{\@tempa\noexpand\thesubeqncnt}}

\makeatother

\allowdisplaybreaks

\begin{document}

\titlepage

\title{Gauss Metric on the Kummer Surface} 
\author{Masahito Hayashi\thanks{masahito.hayashi@oit.ac.jp}\\
Osaka Institute of Technology, Osaka 535-8585, Japan\\
Kazuyasu Shigemoto\thanks{shigemot@tezukayama-u.ac.jp} \\
Tezukayama University, Nara 631-8501, Japan\\
Takuya Tsukioka\thanks{tsukioka@bukkyo-u.ac.jp}\\
Bukkyo University, Kyoto 603-8301, Japan\\
}
\date{\empty}


\maketitle

\abstract{%

On the Kummer surface, we have obtained two different Gauss metrices 
by parametrizing it in two ways.
We have found that these Gauss metrices are not Ricci flat.
The double sphere, which is the special case of the Kummer surface, has the K\"{a}hler 
metric and the first Chern class of it 
does not vanish. Its metric is the Einstein metric which is not Ricci flat.}


\section{Introduction} 
\setcounter{equation}{0}

The research on the Kummer surface has a long history~\cite{Dolgachev}. 
Fresnel found the equation that determines a speed of light travels in biaxial crystal 
in 1822~\cite{Fresnel},  
which was noticed later to be the special case of the Kummer surface.
The equation he discovered can be rewritten as a quartic  equation concerning $(x,y,z)$:
\begin{align}
 & \left( {x}^2 + {y}^2 + {z}^2 \right)\left( {a}^2{x}^2 + {b}^2{y}^2 + {c}^2{z}^2 \right) 
  \notag\\
 &\hspace{2cm}-\left\{
     {a}^2\left( {b}^2 + {c}^2 \right)x^2
    +{b}^2\left( {c}^2 + {a}^2 \right)y^2
    +{c}^2\left( {a}^2 + {b}^2 \right)z^2
   \right\}
  +a^2b^2c^2
  = 0.
  \label{eq:Fresnel}
\end{align}
In 1864, Kummer found an interesting quartic surface, which has a very high 
symmetry and has the interesting structure of the singularities. 
This Kummer surface has 16 nodes and 16 trope-conics and these form 
the interesting structure, which we call $16_6$ structure. Each node lies on 
six trope-conics and each trope-conics contains six nodes~\cite{Kummer}.

For given manifolds/varieties, we can assign several metrices. 
As an example, consider a surface represented by  $x^2+y^2-z^2=-1$ and parametrize 
it in the form
$
  x=\sinh u \cos v, ~y=\sinh u \sin v, ~z=\cosh u.
$
This surface can be considered to be AdS$_2$ 
by identifying $z=t(=\text{time})$. The metric for this surface is
\[
  ds^2=dx^2+dy^2-dt^2=du^2+\sinh^2 u \,dv^2.
\]
On the other hand, this surface can also be considered to be 
a hyperboloid of two sheets. 
The metric in this case is
\[
  ds^2=dx^2+dy^2+dz^2=(\cosh^2u+\sinh^2u)\,du^2+\sinh^2 u \,dv^2.
\]

For the Kummer surface, a similar situation happens.
In the complex differential geometry, the
Kummer surface is the special case of the K3 manifolds.
They are two dimensional complex manifolds
which have K\"{a}hler metric with vanishing both a first Chern 
class $c_1$ and a Hodge number $h^{0,1}$.
The metric of the Kummer surface becomes Ricci flat according 
to the Calabi-Yau theory~\cite{Calabi1,Calabi2,Yau}.  

While, the original Kummer surface, Fresnel's wave surface Eq.(\ref{eq:Fresnel}) is one of the examples, 
is not defined as a complex manifold. Then it is not guaranteed to 
have the K\"{a}hler metric. 
Instead, we have calculated Gauss metric on the original Kummer surface,
which is a non-linear space defined by a quartic relation.
Gauss metric seems to be the systematic metric.

In this paper, we have defined the Gauss metric by two different parametrizations for the 
Kummer surface, and have found that the Gauss metric is not Ricci flat. We also discuss the 
special case, that is, the double sphere and the metric becomes 
the Einstein metric which is not Ricci flat.

\section{
Gauss metric on the Kummer surface} 
\setcounter{equation}{0}
%
\subsection{Kummer surface}
First,  we summarize briefly the Kummer surface coordinated by 
the genus two hyperelliptic $\wp$ function according to 
Baker's book~\cite{Baker1}.

We denote
\begin{equation}
  \wp_{ij} = - \frac{\partial^2}{\partial u_i \partial u_j} \log\sigma(u_1,u_2),\quad
  \wp_{ijk\ell}= \frac{\partial^2}{\partial u_k \partial u_\ell} \wp_{ij}.
\end{equation}
These quantities satisfy the following differential equations~\cite{Baker1,Baker2},  
\begin{align}
&1)\ \wp_{2222} -6\wp_{22}^2 -4\wp_{21} -\lambda_4\wp_{22} - \frac{1}{2}\lambda_3 =0,
\label{2e2}\\
&2)\ \wp_{2221} -6\wp_{22}\wp_{21} + 2\wp_{11} -\lambda_4\wp_{21} =0,
\label{2e3}\\
&3)\ \wp_{2211} -4\wp_{21}^2 -2\wp_{22} \wp_{11} -\frac{1}{2}\lambda_3\wp_{21} =0,
\label{2e4}\\
&4)\ \wp_{2111} -6\wp_{21}\wp_{11} -\lambda_2\wp_{21} +\frac{1}{2}\lambda_1\wp_{22} +\lambda_0 =0,
\label{2e5}\\
&5)\ \wp_{1111} -6\wp_{11}^2 -\lambda_2\wp_{11} -\lambda_1\wp_{21} +3 \lambda_0 \wp_{22}
-\frac{1}{8}\lambda_3 \lambda_1 +\frac{1}{2}\lambda_4 \lambda_0=0  . 
\label{2e6}
\end{align}
Here $\lam_0,\dots,\lam_4$ are coefficients included in a genus two hyperelliptic curve
\begin{equation}
y^2=
f_5(x) =4x^5 + \lam_4x^4 + \lam_3x^3 + \lam_2x^2 + \lam_1x + \lam_0.
\label{2e1}
\end{equation}

We have four consistency conditions
\begin{align*}
   \partial_2\wp_{1111} = \partial_1\wp_{2111},\quad
   \partial_2\wp_{2111} = \partial_1\wp_{2211},\quad
   \partial_2\wp_{2211} = \partial_1\wp_{2221},\quad
   \partial_2\wp_{2221} = \partial_1\wp_{2222}.
\end{align*}
These equations are written as follows
\begin{equation}
\left( \begin{array}{cccc}
-\lam_{0}\hphantom{-}  &  \frac{\lam_1}{2} & 2\wp_{11}  & -2\wp_{21}\hphantom{-} \\
\frac{\lam_1}{2}       & -\lam_2-4\wp_{11} & \hphantom{-}\frac{\lam_3}{2}+2\wp_{21}    & 2\wp_{22}\\
 2\wp_{11}             &  \hphantom{-}\frac{\lam_3}{2}+2\wp_{21} &   -\lam_4-4\wp_{22} & 2\\
-2\wp_{21}\hphantom{-} &  2\wp_{22}        & 2          & 0 \\
\end{array} \right)
\left( \begin{array}{c}
\wp_{222} \\
\wp_{221} \\
\wp_{211} \\
\wp_{111} 
\end{array} \right )=
K
\left( \begin{array}{c}
\wp_{222} \\
\wp_{221} \\
\wp_{211} \\
\wp_{111} 
\end{array} \right )=0 .
\label{2e7}
\end{equation}
A condition in which $\wp_{222}=\partial_2\wp_{22}$, $\wp_{221}=\partial_2\wp_{21}=\partial_1\wp_{22}$, 
$\wp_{211}=\partial_2\wp_{11}=\partial_1\wp_{21}$ and $\wp_{111}=\partial_1\wp_{11}$ 
have a non-trivial solution, is $\det K=0$. 
Rewrite $X=\wp_{22}$, $Y=\wp_{21}$ and $Z=\wp_{11}$.  By taking $\lam_4, \lam_3, \lam_2, \lam_1, \lam_0$ 
to be real, $(X,Y,Z)$ gives the real quartic surface in three dimensional Euclidean space in the form 
\begin{equation}
\det K=
\left|
\begin{array}{cccc}
-\lam_{0}\hphantom{-} &  \frac{\lam_1}{2} & 2Z  & -2Y\hphantom{-} \\
\frac{\lam_1}{2}      & -\lam_2-4X        & \hphantom{-}\frac{\lam_3}{2}+2Y & 2X\\
 2Z                   &  \hphantom{-}\frac{\lam_3}{2}+2Y &   -\lam_4-4X     & 2\\
-2Y\hphantom{-}       &  2X               & 2   & 0 \\
\end{array}
\right| =0 .
\label{eq:Kummer_surface}
\end{equation}
This surface is called the Kummer surface~\cite{Baker2}.

In the following, symbols $u=u_1$ and $v=u_2$ are used. $X$, $Y$ and $Z$ are expressed by
using $\sigma$ function and its derivatives as  
\begin{align}
  X &= \wp_{22}(u,v) = \frac{{\sigma_2}^2-\sigma \sigma_{22}}{\sigma^2},
\label{2e9}\\
  Y &= \wp_{21}(u,v) = \frac{ \sigma_2\sigma_1-\sigma \sigma_{21}}{\sigma^2},
\label{2e10}\\
  Z &= \wp_{11}(u,v) = \frac{{\sigma_1}^2-\sigma \sigma_{11}}{\sigma^2}.
\label{2e11}
\end{align} 
There are various solutions for  Eqs.\eqref{2e2}--\eqref{2e6}~\cite{Hayashi8}. In this note, we 
adopt the following power series expansion of $\sigma$ function with respect to $u$ and $v$ of the form ~\cite{Baker1}
\begin{align}
\sigma(u,v) &= u + \sigma^{(3)}(u,v)+\sigma^{(5)}(u,v)+\sigma^{(7)}(u,v)+\sigma^{(9)}(u,v) +\cdots , 
\label{eq:expansion_sigma}
\end{align}
where
\begin{align}
&\sigma^{(3)}(u,v)=\frac{\lam_2}{24} u^3-\frac{1}{3}v^3, 
\label{2e13}\\
&\sigma^{(5)}(u,v)\notag\\
&=-\frac{1}{5!}
\left(
\left(\frac{\lam_0 \lam_4}{2} -\frac{\lam_1 \lam_3}{8} -\frac{\lam_2^2}{16} \right) u^5 
+10 \lam_0  u^4v
+ 5 \lam_1   u^3v^2
+ 5 \lam_2   u^2v^3
+\frac{5}{2} \lam_3 uv^4 
+2 \lam_4 v^5
\right),
\label{2e14}\\
& \cdots.
\nonumber 
\end{align}
An explicit expression for $\sigma^{(7)}$ is given in Appendix A. 
We must notice that all terms contained in $\sigma^{(i)}\ (i\ge 5)$ include some of
$\lambda_{j}\ (j=0,\dots,4)$.

By this power series solution, the  Kummer surface Eq.\eqref{eq:Kummer_surface} is parametrized by two real parameters $(u,v)$. We have calculated $\det K$ order by order of $\sigma$ in Eq.\eqref{eq:expansion_sigma}. As a result, the following {\color{blue}} were found
\begin{align}
&{\rm\hphantom{ii}i}) ~\text{For}\ \sigma=u+\sigma^{(3)}, ~\sigma^8 \det K = 0 ~\text{up to 5-th order}, 
\label{2e15}\\
&{\rm\hphantom{i}ii}) ~\text{For}\ \sigma=u+\sigma^{(3)}+\sigma^{(5)}, ~\sigma^8 \det K =0 ~\text{up to 7-th order}, 
\label{2e16}\\
&{\rm iii}) ~\text{For}\ \sigma=u+\sigma^{(3)}+\sigma^{(5)}+\sigma^{(7)}, ~\sigma^8 \det K =0 ~\text{up to 9-th order}.
\label{2e17}
\end{align}

\subsection{Gauss metric on the Kummer surface}
By using the method of Gauss, we define the metric on 
the Kummer surface through the first fundamental form~\cite{Klingenberg}.
We denote the Kummer surface as ${\bf S}(u,v)=(X(u,v), Y(u,v), Z(u,v))$, and define
\begin{align}
\partial_1 {\bf S}=\partial_u {\bf S}&=(\partial_u X, \partial_u Y, \partial_u Z)
=(\wp_{221},\wp_{211}, \wp_{111}),
\label{2e18}\\
\partial_2 {\bf S}=\partial_v {\bf S}&=(\partial_v X, \partial_v Y, \partial_v Z)
=(\wp_{222},\wp_{221}, \wp_{211}).
\label{2e19}
\end{align}
Here $\wp_{ijk}$'s are expressed as
\begin{align}
&\wp_{222}=-\frac{1}{\sigma^3}\left( \sigma^2 \sigma_{222} -3 \sigma \sigma_2 \sigma_{22}
+2 {\sigma_2}^3\right),
\label{2e24}\\
&\wp_{221}=-\frac{1}{\sigma^3} \left( \sigma^2 \sigma_{221} 
-\sigma (\sigma_{22} \sigma_{1}+2\sigma_{21} \sigma_{2})+2{\sigma_2}^2 \sigma_1\right),
\label{2e25}\\
&\wp_{211}=-\frac{1}{\sigma^3}\left( \sigma^2 \sigma_{211} 
-\sigma (2 \sigma_{21} \sigma_{1}+\sigma_{11} \sigma_{2})+2\sigma_2 {\sigma_1}^2\right),
\label{2e26}\\
&\wp_{111}=-\frac{1}{\sigma^3}\left( \sigma^2 \sigma_{111} -3 \sigma \sigma_1 \sigma_{11}
+2 {\sigma_1}^3\right).
\label{2e27}
\end{align}
Then Gauss metric is defined as  
\begin{align}
  ds^2 &= \partial_u {\bf S} \cdot \partial_u {\bf S}\, du^2
        +2\partial_u {\bf S} \cdot \partial_v {\bf S}\, du dv
         +\partial_v {\bf S} \cdot \partial_v {\bf S}\, dv^2
\nonumber\\
      &=g_{11}\,du^2 + 2g_{12}\,du dv + g_{22}\,dv^2,
\label{2e20}
\end{align}
which gives 
\begin{align}
&g_{11}=\wp_{221}^2+\wp_{211}^2+\wp_{111}^2,
\label{2e21}\\
&g_{12}=g_{21}=\wp_{222}\wp_{221}+\wp_{221}\wp_{211}+\wp_{211}\wp_{111},
\label{2e22}\\
&g_{22}=\wp_{222}^2+\wp_{221}^2+\wp_{211}^2.
\label{2e23}
\end{align}
By using the power series solution, $g_{ij}$'s and $\det g$ are given as
%
%
\begin{align*}
  g_{11} 
  =& \frac{4}{\sigma^6} \left( 1 + u^2v^2 + v^4 + \frac{4}{3}uv^5+\frac{4}{9} v^8 + \text{$\lam_{i}$ dependent terms} \right)
  \equiv\frac{\widehat{g}_{11}}{\sigma^6},
\\
  g_{12} 
  =& \frac{-4}{\sigma^6} \left( v^2 + u^3v + uv^3 + 3u^2v^4 + \frac{2}{3}v^6 + \frac{5}{3} u v^7 + \frac{2}{27}v^{10}+\text{$\lam_{i}$ dependent terms} \right)
  \equiv\frac{\widehat{g}_{12}}{\sigma^6},
\\
  g_{22} 
  =& \frac{4}{\sigma^6} \left( u^4 + u^2v^2 + v^4 + \frac{14}{3}u^3v^3 + \frac{4}{3}uv^5 + \frac{17}{3}u^2 v^6 + \frac{4}{9}v^8 + \frac{14}{27}uv^9 + \frac{1}{81} v^{12} \right.
\\
&\hphantom{\frac{4}{\sigma^{12}}\left(\right.}+\text{$\lam_{i}$ dependent terms}\left.\vphantom{\frac12} \right)
  \equiv\frac{\widehat{g}_{22}}{\sigma^6},
\end{align*}
\begin{align*}
 \det g 
  =& \frac{16}{\sigma^{12}} \left( u^4 + u^2v^2 + \frac{8}{3}u^3v^3 - \frac{2}{3}uv^5 + \frac{2}{3}u^2 v^6 + \frac{1}{9}v^8 - \frac{40}{27}u v^9 + u^2 v^{10} + \frac{25}{81}v^{12} \right.
\nonumber\\
&\hphantom{\frac{16}{\sigma^{12}}\left(\right.}- \frac{2}{3}u v^{13} + \frac{1}{9} v^{16} + \text{$\lam_{i}$ dependent terms} \left.\vphantom{\frac12} \right)
  \equiv\frac{\widehat{D}}{\sigma^{12}} ,
\end{align*}
and $g^{ij}$'s are found as follows
\begin{align}
&g^{11}=\frac{g_{22}}{\det g}=\widehat{g}_{22}\frac{\sigma^6}{\widehat{D}},
\label{2e32}\\
&g^{12}=g^{21}=-\frac{g_{12}}{\det g}=-\widehat{g}_{12}\frac{\sigma^6}{\widehat{D}},
\label{2e33}\\
&g^{22}=\frac{g_{11}}{\det g}=\widehat{g}_{11}\frac{\sigma^6}{\widehat{D}}.
\label{2e34}
\end{align}

In the standard way, we calculate 
\begin{align}
&{\MGamma^{\lam}}_{\mu \nu}=\frac{1}{2} g^{\lam \sigma}
( \partial_{\mu} g_{\sigma \nu}+ \partial_{\nu} g_{\mu \sigma}
- \partial_{\sigma} g_{\mu \nu}),
\label{2e35}\\
&{R^{\alpha}}_{\beta \mu \nu}=
  \partial_{\mu} {\MGamma^{\alpha}}_{\beta \nu}
- \partial_{\nu} {\MGamma^{\alpha}}_{\beta \mu}
+{\MGamma^{\alpha}}_{\tau \mu} {\MGamma^{\tau}}_{\beta \nu}
-{\MGamma^{\alpha}}_{\tau \nu} {\MGamma^{\tau}}_{\beta \mu},
\label{2e36}\\
&R_{\mu \nu}={R^{\alpha}}_{\mu \alpha \nu}.
\label{2e37}
\end{align}
We denote components of Ricci tensor in the form  
\begin{align}
R_{11}=\frac{\widehat{R}_{11}}{\sigma^2 {\widehat{D}}^2}, \quad
R_{12}=\frac{\widehat{R}_{12}}{\sigma^2 {\widehat{D}}^2}, \quad
R_{22}=\frac{\widehat{R}_{22}}{\sigma^2 {\widehat{D}}^2}.
\label{2e41}
\end{align}
Then $\widehat{R}_{ij}$'s become infinite power series. The lowest order of $\widehat{R}_{11}$, 
$\widehat{R}_{12}$ and $\widehat{R}_{22}$ are 10, 12 and 14, respectively. 

The $\lam_i (i=0,\dots,4)$ independent lowest terms of 
$\widehat{R}_{11}$, $\widehat{R}_{12}$ and $\widehat{R}_{22}$ are
\begin{align}
\text{$\lam_i$ independent lowest terms of }
\widehat{R}_{11}&=-2^{10}\, u^5 v^5,
\label{2e42}\\
\text{$\lam_i$ independent lowest terms of }
\widehat{R}_{12}&=\hphantom{-}2^{10}\, u^5 v^{7},
\label{2e43}\\
\text{$\lam_i$ independent lowest terms of }
\widehat{R}_{22}&=-2^{10}\, u^5 v^5 
\left( u^4+u^2 v^{2}+ v^{4} \right).
\label{2e44}
\end{align}
We cannot eliminate these terms by adding higher $\lam_i (i=0,\dots,4)$ dependent terms. 
Indeed, we have checked that these terms are unchanged even if we take
i) $\sigma=u+\sigma^{(3)}$, ii) $\sigma=u+\sigma^{(3)}+\sigma^{(5)}$,  iii) $\sigma=u+\sigma^{(3)}+\sigma^{(5)}+\sigma^{(7)}$.

\subsection{Another Gauss metric on the Kummer surface}
 In connection with the Jacobi's inversion problem, 
we have another parametrization
of the Kummer surface\cite{Baker3}.

By using the relations, 
$$
\int^{x_1}\!\! \frac{  dx}{\sqrt{f_5(x)}}+\int^{x_2}\!\! \frac{  dx}{\sqrt{f_5(x)}}=u, \quad
\int^{x_1}\!\! \frac{x dx}{\sqrt{f_5(x)}}+\int^{x_2}\!\! \frac{x dx}{\sqrt{f_5(x)}}=v,
$$
where $f_5(x)$ is given in Eq.\eqref{2e1}, 
we express the symmetric combination of $x_1, x_2$ as the function of $u, v$, 
which is the inversion problem of expressing $u, v$ as the function of $x_1, x_2$.
Then genus two hyperelliptic $\wp$ functions are expressed as the symmetric 
function of $x_1, x_2$. 

Thus we obtain 
\begin{align}
&X=\wp_{22}=x_1+x_2,\quad
 Y=\wp_{21}=-x_1x_2,\quad
 Z=\wp_{11}=\frac{F(x_1,x_2)-2 y_1 y_2}{4(x_1-x_2)^2},
\label{2e45}\\
&F(x_1,x_2)=4x_1^2x_2^2(x_1+x_2) + 2\lam_{4}x_1^2 x_2^2 + \lam_{3}x_1x_2(x_1+x_2)
+ 2\lam_{2}x_1x_2 + \lam_{1}(x_1+x_2) + 2\lam_{0}, 
\label{2e46}\\
&y_i^2=f_5(x_i)=4 x_i^5+\lam_4 x_i^4+\lam_3 x_i^3+\lam_2 x_i^2+\lam_1 x_i+\lam_0,\ (i=1, 2),
\label{2e47}
\end{align}
We have checked that these $X, Y$ and $Z$ satisfy Eq.\eqref{eq:Kummer_surface}.

We define the vector
\[
 {\bf S}(x_1,x_2)=(X(x_1,x_2), Y(x_1,x_2), Z(x_1,x_2))=(x_1+x_2, -x_1 x_2, Z(x_1,x_2)),
\]
and
\begin{align}
\frac{\partial {\bf S}(x_1,x_2)}{\partial x_1}
=\left(1,-x_2, \frac{\partial Z(x_1,x_2)}{\partial x_1}\right), \quad
\frac{\partial {\bf S}(x_1,x_2)}{\partial x_2}
=\left(1,-x_1, \frac{\partial Z(x_1,x_2)}{\partial x_2}\right).
\label{2e48}
\end{align}
Then Gauss metric is defined as  
\begin{align}
ds^2&=
\frac{\partial {\bf S}}{\partial x_1} \cdot \frac{\partial {\bf S}}{\partial x_1}\, {dx_1}^2
+2\,\frac{\partial {\bf S}}{\partial x_1} \cdot \frac{\partial {\bf S}}{\partial x_2} \, dx_1 dx_2
+\frac{\partial {\bf S}}{\partial x_2} \cdot \frac{\partial {\bf S}}{\partial x_2}\, {dx_2}^2
\nonumber\\
&=g_{11}(x_1,x_2)\, {dx_1}^2+2g_{12}(x_1,x_2)\, dx_1 dx_2+g_{22}(x_1,x_2)\, {dx_2}^2,
\label{2e49}
\end{align}
which gives 
\begin{align}
&g_{11}=1+x_2^2+\left({\frac{\partial Z}{\partial x_1}}\right)^2,
\label{2e50}\\
&g_{12}=g_{21}=
1+x_1 x_2+{\frac{\partial Z}{\partial x_1}}{\frac{\partial Z}{\partial x_2}},
\label{2e51}\\
&g_{22}=1+x_1^2+\left({\frac{\partial Z}{\partial x_2}}\right)^2.
\label{2e52}
\end{align}
Here
\begin{align}
\frac{\partial Z}{\partial x_1}
&=
-\frac{F-2 y_1 y_2}{2(x_1-x_2)^3}
+
\frac{1}{4 (x_1-x_2)^2}
 \left(
   \frac{\partial F}{\partial x_1} - \frac{y_2}{y_1}\left(\frac{d y_1^2}{d x_1} \right)
 \right),
\\
\frac{\partial Z}{\partial x_2}
&=\hphantom{-}
\frac{F-2 y_1 y_2}{2(x_1-x_2)^3}
+
\frac{1}{4 (x_1-x_2)^2}
 \left(
   \frac{\partial F}{\partial x_2} - \frac{y_1}{y_2}\left(\frac{d y_2^2}{d x_2}\right)
 \right).
\end{align}

We have calculated Ricci tensors $R_{11}$, $R_{12}$ and $R_{22}$. 
They are not equal to zero. The formulas are very complex and we avoid showing them in detail here.
Instead, we show the special case where $\lam_0 \ne 0, ~\lam_1=\lam_2=\lam_3=\lam_4=0$.
\[
  R_{11}=N_{11}/D_{11},\quad R_{22}=N_{22}/D_{22},\quad R_{12}=N_{12}/D_{12}, 
\]
\begin{align}
&N_{11}
=\sqrt{4 x_2^5+\lam_0} \,(-2^{20} \cdot 3 \cdot 5  x_1^{30} x_2^{46}+\cdots)
+\sqrt{4 x_1^5+\lam_0} \,( 2^{20} x_1^{27} x_2^{49}+\cdots) ,
\nonumber\\
&D_{11}
=\sqrt{4 x_2^5+\lam_0} \,( 2^{18} x_1^{33} x_2^{47}+\cdots)
+\sqrt{4 x_1^5+\lam_0} \,(-2^{22} x_1^{31} x_2^{49}+\cdots) ,
\nonumber\\
&N_{22}
=\sqrt{4 x_2^5+\lam_0} \,(-2^{20} \cdot 3 x_1^{27} x_2^{49}+\cdots)
+\sqrt{4 x_1^5+\lam_0} \,( 2^{18} x_1^{24} x_2^{52}+\cdots) ,
\nonumber\\
&D_{22}
=\sqrt{4 x_2^5+\lam_0} \,( 2^{18} x_1^{28} x_2^{52}+\cdots)
+\sqrt{4 x_1^5+\lam_0} \,(-2^{22} x_1^{26} x_2^{52}+\cdots) ,
\nonumber\\
&N_{12}
=\sqrt{4 x_1^5+\lam_0} \sqrt{4 x_2^5+\lam_0} \,(-2^{18} \cdot 3^3 x_1^{26} x_2^{50}+\cdots)
+2^{21} x_1^{28} x_2^{53}+\cdots ,
\nonumber\\
&D_{12}
=\sqrt{4 x_1^5+\lam_0} \sqrt{4 x_2^5+\lam_0} \,( 2^{18} x_1^{28} x_2^{52}+\cdots)
-2^{24} x_1^{31} x_2^{54}+\cdots .
\end{align}

\subsection{Special Kummer surface : double sphere}

There is another way to represent the Kummer surface, which is the identity relation for the 
genus two hyperelliptic $\vartheta$ functions in the form~\cite{Gopel,Hudson}
\begin{align}
X^4&+Y^4+Z^4+T^4+A\left(X^2 T^2+Y^2 Z^2\right)+B\left(Y^2 T^2+Z^2 X^2\right)+C\left(Z^2 T^2+X^2 Y^2\right)
\nonumber\\
&+2 D X Y Z T=0,
\label{2e55}\\
&A=\frac{\beta^4 +\gamma^4-\alpha^4-\delta^4}{\alpha^2 \delta^2-\beta^2  \gamma^2}, \quad
 B=\frac{\gamma^4+\alpha^4-\beta^4 -\delta^4}{\beta^2  \delta^2-\gamma^2 \alpha^2}, \quad
 C=\frac{\alpha^4+\beta^4- \gamma^4-\delta^4}{\gamma^2 \delta^2-\alpha^2 \beta^2 },
\notag\\
&D=\frac{ \alpha \beta \gamma \delta (2-A)(2-B)(2-C)}{(\alpha^2+\beta^2+\gamma^2+\delta^2)^2}.
\notag
\end{align}
$(X,Y,Z,T)$ is a homogenious coordinate of $\mathbb{P}^3$.
If $\alpha=\beta=\gamma=1$, $\delta=\sqrt{-3}$, then $A=B=C=2$, $D=0$ which gives 
\begin{equation}
\left(X^2+Y^2+Z^2+T^2\right)^2=0.
\label{2e58}
\end{equation}
By rewriting $x=i X/T, \,y=i Y/T, \,z=i Z/T$, we obtain the double sphere as~\cite{Endrass}
\begin{equation}
\left(x^2+y^2+z^2-1\right)^2=0.
\label{2e59}
\end{equation}
Eq.\eqref{2e59} is also derived from Fresnel's wave surface Eq.\eqref{eq:Fresnel} by placing $a=b=c=1$.

The metric of this Kummer surface is nothing but the one of sphere
$x^2+y^2+z^2-1=0$. By using the spherical coordinate $x=\sin \theta \cos \phi,\, 
y=\sin \theta \sin \phi,\, z=\cos \theta$, the metric of the sphere is given by 
\begin{equation}
ds^2=d \theta^2+\sin^2 \theta\ d \phi^2,
\label{2e60}
\end{equation}
which gives the Einstein metric $R_{ij}=g_{ij}$ with the constant  scalar 
curvature $R=2$. Just as in the Riemann sphere, 
we can introduce the complex coordinate in the form 
$\zeta=(x+i y)/(1-z)=\cot( \theta/2) e^{i \phi}=u+i v$.
By using the K\"{a}hler potential $\mathit\Phi=4 \log (1+\zeta \bar{\zeta})$, 
we can rewrite the metric to the K\"{a}hler metric form 
\begin{equation}
ds^2=\frac{\partial^2 \mathit\Phi }{\partial \zeta \partial \bar{\zeta} }
d\zeta d\bar{\zeta}
=\frac{4 d\zeta d\bar{\zeta}} {(1+ \zeta \bar{\zeta})^2}=
\frac{4 (du^2+dv^2) } {(1+ u^2+v^2)^2}\,.
\label{2e61}
\end{equation} 
Then the first Chern class $c_1$ is calculated in the form
\begin{equation}
c_1=\frac{i}{\pi} \int \frac{d \zeta \wedge d \bar{\zeta}}
{(1+ \zeta \bar{\zeta})^2}=2.
\label{2e62}
\end{equation} 
%

\section{Summary and Discussions} 
\setcounter{equation}{0}


We investigated the metric of Kummer surface Eq.(\ref{eq:Kummer_surface}) coordinated 
by the genus two hyperelliptic $\wp$ functions in two ways.
In the first method, we used the expression of the genus two hyperelliptic $\wp$ 
functions with $u$ and $v$ as variables by solving differential equations Eqs(2.2)--(2.6), where 
the solution was represented as the power series of the $\sigma$ function.
In the second method, by using the idea of the Jacobi's inversion problem, 
we used the expression of  genus two hyperelliptic $\wp$ functions
as the symmetric function of $x_1$ and $x_2$.
We found that these Gauss metrices are not Ricci flat.

Next, we considered the metric on the double sphere, which is the special case of Kummer surface.
By using the parametrization of the Riemann sphere, the metric of this double sphere 
can be written in the K\"{a}hler metric. The first Chern class of this  K\"{a}hler metric
does not vanish. This metric becomes the Einstein metric which is not Ricci flat.

\begin{appendices}
\setcounter{equation}{0}
\section{ The explicit expression of $\sigma^{(7)}(u,v)$} 
\begin{align}
\sigma^{(7)}(u,v) &= \frac{1}{7!}\sum_{k=0}^7 \binom{7}{k} \, h_k u^kv^{7-k},
\end{align}
where
\begin{align}
  h_0 &=-\lam_3-2 \lam_4^2,
\label{A2}\\
  h_1 &=-2 \lam_2-\frac{1}{2} \lam_3 \lam_4,
\label{A3}\\
  h_2 &=-2 \lam_1-\frac{1}{2} \lam_2 \lam_4,
\label{A4}\\
  h_3 &=-2 \lam_0-\frac{1}{8} \lam_2 \lam_3-\frac{1}{2} \lam_1 \lam_4,
\label{A5}\\
  h_4 &=-\frac{1}{4} \lam_1 \lam_3-\lam_0 \lam_4 -\frac{1}{8} {\lam_2}^2,
\label{A6}\\
  h_5 &=-\frac{3}{2} \lam_0 \lam_3-\frac{1}{4} \lam_1 \lam_2,
\label{A7}\\
  h_6 &=-\frac{11}{2} \lam_0 \lam_2+ {\lam_1}^2,
\label{A8}\\
  h_7 &=\frac{1}{64} {\lam_2}^3+\frac{3}{32} \lam_1 \lam_2 \lam_3-\frac{15}{8} \lam_0 \lam_2 \lam_4-\frac{1}{2} \lam_0 \lam_1+\frac{3}{8} \lam_0 {\lam_3}^2+\frac{3}{8} {\lam_1}^2 \lam_4,
\label{A9}
\end{align}
%
\end{appendices}


\end{document}